\documentclass[10pt]{elsarticle}
\usepackage{caption}
\usepackage{subcaption}
\usepackage{float}
\usepackage{booktabs}  
\usepackage{bbding}
\usepackage{graphicx}
\usepackage{xcolor}%
\usepackage{physics}
\usepackage{booktabs}
\graphicspath{ {./charts/} }
\usepackage{amsthm,amsmath,amssymb}
\usepackage[linesnumbered,ruled]{algorithm2e}

\begin{document}
\begin{frontmatter}
\title{Weak Random Feature Method for Solving Partial Differential Equations}
\author[1]{Mikhail Kuvakin}
\ead{kuvakin.ms@phystech.edu}

\author[2]{Zijian Mei\corref{cor1}}
\ead{meizijian518@ustc.mail.edu.cn}

\author[3]{Jingrun Chen}
 \ead{jingrunchen@ustc.edu.cn}

\address[1]{Moscow Institute of Physics and Technology}
\address[2]{School of Artificial Intelligence and Data Science and Suzhou Institute
for Advanced Research, University of Science and Technology of China, Suzhou, China}
\address[3]{School of Mathematical Sciences and Suzhou Institute for Advanced
Research, University of Science and Technology of China, and Suzhou Big Data \& AI Research and Engineering Center, Suzhou, China}

\cortext[cor1]{Corresponding author}

\begin{abstract}
The random feature method (RFM) has demonstrated great potential in bridging traditional numerical methods and machine learning techniques for solving partial differential equations (PDEs). It retains the advantages of mesh-free approaches while achieving spectral accuracy for smooth solutions, without the need for iterative procedures. However, the implementation of RFM in the identification of weak solutions remains a subject of limited comprehension, despite crucial role of weak solutions in addressing numerous applied problems. While the direct application of RFM to problems without strong solutions is fraught with potential challenges, we propose an enhancement to the original random feature method that is specifically suited for finding weak solutions and is termed as Weak RFM. Essentially, Weak RFM reformulates the original RFM by adopting the weak form of the governing equations and constructing a new linear system through the use of carefully designed test functions, ensuring that the resulting solution satisfies the weak form by default. To rigorously evaluate the performance of the proposed method, we conduct extensive experiments on a variety of benchmark problems, including challenging three-dimensional cases, and compare its performance with state of the art machine learning-based approaches. The results demonstrate that Weak RFM achieves comparable or superior accuracy while significantly reducing computational time and memory consumption, highlighting its potential as a highly efficient and robust tool for finding weak solutions to various PDE problems.
\end{abstract}
\end{frontmatter}

\section{Introduction}
\label{sec:intro}

Most conceivable systems can be described by a set of parameters that are interrelated through equations, often differential ones. Solving these equations enables the identification of explicit relationships between parameters, allowing prediction of behavior of one parameter based on changes in others. This straightforward paradigm is pervasive across the natural sciences. In physics, for instance, differential equations can describe gas density in a core of star \cite{chandrasekhar}, oscillations in chemical concentrations \cite{zhabotinsky}, and the rate of infection spread in biology \cite{hethcote}. Even in fields beyond natural sciences, such as economics, this approach is applied to evaluate the value of options and other securities \cite{black}. For this reason, solving differential equations is a cornerstone of applied mathematics. While some solutions can be derived analytically and others, such as the Schrödinger equation, allow for asymptotic solutions \cite{maslov}, the majority of problems lack explicit solution formulas.

For such cases, numerous numerical methods have been developed, broadly categorized into traditional and machine learning (ML)-based approaches. Traditional methods, such as the finite difference method (FDM) and finite element method (FEM), have a long history that gained momentum with the development of early computers. These methods simplify the process of solving differential equations by breaking it down into iterative steps or solving a linear system. Although well-established in both theory and practice, traditional methods often face limitations with complex domains and high-dimensional or high-order problems.

Meanwhile, along with the development of machine learning, the idea of using neural networks for solving equations has evolved \cite{lee1990, meade1994, lagaris1998}. The concept remained in a relatively nascent state until the combination of advances in computing power and machine learning led to almost simultaneous emergence of several comparable ML-based methods \cite{deepritz, deepgal, pinn}. These methods represent the solution as a neural network, transforming the task of solving a differential equation into an optimization problem addressable by various optimizers. Such methods are readily applicable to complex domains and have a potential to avoid the curse of dimensionality. Furthermore, continuous advancements in deep learning have led to various enhancement strategies, such as loss balancing \cite{relobralo}, transfer learning \cite{transfer}, and novel optimizers \cite{nncg}. Nonetheless, these methods are not yet capable of fully replacing traditional approaches, as they often fall short of meeting industrial demands for accuracy, computational efficiency or stability.

The random feature method (RFM) \cite{chen2022bridging} bridges the gap between traditional numerical approaches and machine learning techniques, effectively combining the strengths of both. From one perspective, RFM resembles ML-based methods by representing the solution as a multi-layer perceptron (MLP) with a single hidden layer, providing a continuous function that eliminates the need for interpolation and easily accommodates complex geometries. From another perspective, RFM bypasses the time-consuming neural network optimization process by solving a linear system via least squares, thereby achieving computational efficiency and accuracy comparable to traditional numerical methods. This dual approach makes RFM a compelling option for efficiently solving differential equations. Various enhancements have been proposed to expand the capabilities of RFM in solving PDEs, adapting it to handle time-dependent \cite{chen2023time} and interface \cite{chen2024interface, mei2024solving} problems. Additionally, new strategies for solving least-squares problems have been introduced to improve the performance of RFM in complex geometries \cite{chen2024methods}.

Given the potential of RFM in finding classical solutions \cite{chen2024optimization}, adapting it for problems lacking solutions in the usual sense appears promising.  Such issues are common in physics \cite{kawashima1999}, particularly in fluid mechanics \cite{spalart1986, lundmark2007}, and often stem from boundary conditions that may be redundant, incompatible, or insufficiently smooth. Notably, the absence of a classical solution does not necessarily indicate that a solution is altogether lacking; rather, a solution that satisfies the equation almost everywhere may still exist. Solutions of this nature are known as weak solutions, and finding them requires applying a weak formulation of the original equation.

However, because the basis functions in RFM are inherently smooth, directly applying RFM to problems with low regularity leads to poor results and introduces a range of challenges. Therefore, this work aims to enhance the RFM algorithm to effectively handle weak solutions. To this end, we modify the original framework and, for simplicity, refer to the proposed method as WRFM. The core idea of Weak RFM is closely aligned with the RFM framework: the solution is represented as an MLP with one hidden layer, and, in contrast to the training process employed in other ML-based methods, a linear system is solved using the least-squares method. The solution to the linear system yields the coefficients of the MLP output level, thereby establishing the solution to the equation. Nevertheless, due to the adoption of the weak formulation of the governing equations, the resulting linear system in WRFM differs significantly from that of RFM. Given that a weak formulation necessitates the satisfaction of an equation for an infinite set of test functions, this could potentially result in a challenge during the implementation of the method. However, drawing inspiration from the approximation theorem \cite{weierstrass1885}, we circumvent this issue by employing a finite set of predefined functions to approximate the test functions.

The subsequent sections of this paper are structured as follows. First, we remind the main ideas of the standard RFM in Section \ref{sec:method}. Next, in Section \ref{sec:new} we introduce Weak RFM and describe in detail its features and differences from the original method. In Section \ref{sec:num_results} we provide several numerical experiments in 2D and 3D cases, compare Weak RFM with existing methods and evaluate how hyperparameters can affect accuracy. Finally, Section \ref{sec:conclusion} inevitably summarizes the results of this work.

\section{Standard RFM}
\label{sec:method}
In this section, we provide a concise overview of the standard RFM, focusing on its partition of unity, loss function, and optimization approach. Let us consider the same problem as for original RFM method:
\begin{equation}\label{eq:1}
\left\{
\begin{aligned}
\mathcal{L} u(\boldsymbol{x})&=f(\boldsymbol{x}) \quad \boldsymbol{x}\in\Omega\,,\\
\mathcal{B} u(\boldsymbol{x})&=g(\boldsymbol{x}) \quad \boldsymbol{x}\in\partial\Omega\,,\\
\end{aligned}
\right.
\end{equation}
where $\boldsymbol{x}\in\mathbb{R}^{d}$, $\Omega\subset\mathbb{R}^{d}$ is bounded domain, $f\colon \Omega \to \mathbb{R}$, $g\colon \partial\Omega \to \mathbb{R}$ and $\mathcal{L}$, $\mathcal{B}$ are linear differential operators:
\begin{equation*}
\mathcal{L} = \sum_{\boldsymbol{\alpha}\in \mathbb{N}^{d}} a_{\boldsymbol{\alpha}}(\boldsymbol{x})\cdot \partial^{\boldsymbol{\alpha}},\quad
\mathcal{B} = \sum_{\boldsymbol{\beta}\in \mathbb{N}^{d}} b_{\boldsymbol{\beta}}(\boldsymbol{x})\cdot \partial^{\boldsymbol{\beta}}.
\end{equation*}

\subsection{Form of Solution}
\label{subsec:form}
Following the approach, introduced in \cite{chen2022bridging}, local solution of equation is being approximated by linear combination of M random feature functions $\{\phi_m\}$:
\begin{equation*}
u_M(\boldsymbol{x})=\sum_{m=1}^{M} u_{m} \phi_{m}(\boldsymbol{x})\,.
\end{equation*}
Each random feature function is constructed by
\begin{equation}\label{phi_m}
\phi_{m}(\boldsymbol{x})=\sigma\left(\boldsymbol{k}_{m} \cdot \boldsymbol{x}+b_{m}\right)\,,
\end{equation}
where $\sigma$ is nonlinear activation, generally one of trigonometric or hyperbolic functions. Elements $\{\boldsymbol{k}_{m}\}$ and $\{b_{m}\}$ are uniformly distributed in the interval $[-R,R]$ and remain fixed until the end of the algorithm, hence the only changeable parameters are $\{u_m\}$.\\
As a result, there are only $M$ degrees of freedom for the approximation space. To handle this feature and incorporate all minutiae of global solution, RFM employs the partition of unity.

\subsection{Partition of Unity}
\label{subsec:partition}
Main idea of this approach is that initial domain $\Omega$ is decomposed into $S$ subdomains and on each subdomain separate solution is being found in normalized coordinates. Local solutions are combined into a single global solution using partition of unity functions.\\
Specifically, a single center point is selected for each subdomain $\left\{\boldsymbol{x}_n\right\}_{n=1}^{S}$ and coordinates are normalized using basic linear transforms:
\begin{equation*}
\begin{aligned}
&l_{ni}(x_i)=\frac{\left(x_i-x_{ni}\right)}{r_{ni}}\,, \quad n=1 \cdots S\,, \quad i=1 \cdots d\,,\\
&\boldsymbol{l}_n(\boldsymbol{x})= (l_{n1}(x_{1}),l_{n2}(x_{2}), \cdots, l_{nd}(x_{d}))^{T}\,,
\end{aligned}
\end{equation*}
where $r_{ni}$ denotes size of $n$-th subdomain in $i$-th dimension. Therefore  $\boldsymbol{l}_n(\boldsymbol{x})$ maps $n$-th subdomain onto $[-1,1]^d$ and allows us to use same partition of unity functions in all subdomains. There are three commonly used types of these functions in one-dimensional case:
\begin{align*}
&\psi^a(x)=  \mathbb{I}_{[-1,1]} ( x )\,, \\
&\psi^b(x)=  \mathbb{I}_{\left[-\frac{5}{4},-\frac{3}{4}\right]} \left( x \right)\cdot \frac{1 + \sin \left(2 \pi x \right)}{2}+\mathbb{I}_{\left[-\frac{3}{4}, \frac{3}{4}\right]}\left( x \right) + \mathbb{I}_{\left[\frac{3}{4}, \frac{5}{4}\right]}\left( x \right)\cdot \frac{1 - \sin \left(2 \pi x \right)}{2}\,, \\
&\psi^c(x)=  \frac{1}{1+e^{-\frac{x+1}{\alpha}}}\cdot\frac{1}{1+e^{-\frac{1-x}{\alpha}}}\,, \\
\end{align*}
where $\mathbb{I}$ is the signal function. Here $\psi^a$ is discontinuous,  while $\psi^b \in C^{1} (\mathbb{R})$ and $\psi^c \in C^{\infty} (\mathbb{R})$. For definition in higher dimensions, tensor product of one-dimensional functions can be used.Thus global solution is being approximated by \eqref{u_M}:
\begin{equation}\label{u_M}
u_M(\boldsymbol{x})=\sum_{n=1}^{S} \psi(\boldsymbol{l}_n(\boldsymbol{x})) \sum_{j=1}^{J_n} u_{n j} \phi_{n j}(\boldsymbol{l}_n(\boldsymbol{x}))\,,
\end{equation}
where $J_n$ is amount of random feature functions in $n$-th subdomain. Total number of degrees of freedom is $M=\sum_{n=1}^{S} J_n$.

\subsection{Loss Function and Optimization for Strong Formulation}
\label{subsec:classic_rfm}
This section describes loss function and consequent optimization process for the case of classic RFM which employs strong formulation of equation. Let us consider two sets of collocation points: $C_I$ and $C_B$, representing correspondingly internal and boundary points of $\Omega$. Furthermore, let us define $E_I$ as the number of governing equations and $E_B$ as the number of boundary conditions. As the strong form of equation is under consideration, loss function have the following form:
\begin{equation*}
\begin{aligned}
    \text { Loss }=\sum_{\boldsymbol{x}_i \in C_I} \sum_{k=1}^{E_I} \left\| \mathcal{L}^k u(\boldsymbol{x}_i)-f^k(\boldsymbol{x}_i) \right\|_{l^2}^2 \\
    +\sum_{\boldsymbol{x}_j \in C_B} \sum_{\ell=1}^{E_B} \left\|\mathcal{B} u(\boldsymbol{x_j})-g^l(\boldsymbol{x_j}) \right\|_{l^2}^2\,.
    \end{aligned}
\end{equation*}
And total number of conditions imposed on $u_M$ is $N=E_I \cdot C_I + E_B \cdot C_B$. Due to the solution design, standard algorithms for least-squares approximation can be employed to solve this optimization problem. One of the most effective way to find parameters $\left\{u_{nj}\right\}_{j=1}^{J_n}$ for all $S$ blocks is to compose a linear system, which incorporates governing equations as well as boundary conditions simultaneously. In case $d=1$, $E_I=1$, $E_B=1$ (for simplicity), the system has the form \eqref{eq:2}:
\begin{equation}\label{eq:2}
A \times \boldsymbol{u}=\boldsymbol{v}\,,
\end{equation}
where 
\begin{equation*}
\begin{aligned}
A=\begin{pmatrix}
A_{1}\\
A_{2}\\
A_{3}\\
A_{4}
\end{pmatrix},&\quad
\boldsymbol{u}=\begin{pmatrix}
u_{11}\\
u_{21}\\
\vdots \\
u_{J_{S} S}
\end{pmatrix},\quad
\boldsymbol{v}=\begin{pmatrix}
\boldsymbol{v_{1}} \\
\boldsymbol{v_{2}} \\
\boldsymbol{v_{3}} \\
\boldsymbol{v_{4}}
\end{pmatrix}\,.
\end{aligned}
\end{equation*}
Let us denote amount of points in $n$-th domain as $P_n$ and total amount of points as $P=\sum_{n=1}^{S} P_n$. Then structure of $A$ and $\boldsymbol{v}$ takes the following form:
\begin{equation*}
\begin{aligned}
\\
&A_1=\begin{pmatrix}
(\mathcal{L}\phi_{1 j}(x_i^{1}))_{i=1,\;j=1}^{P_1,\;J_1} & 0 & \ldots & 0\\
0 & (\mathcal{L}\phi_{2 j}(x_i^{2}))_{i=1,\;j=1}^{P_2,\;J_2} &\ldots & 0\\
\vdots & \vdots & \ddots & \vdots\\
0 &0 & \ldots & (\mathcal{L}\phi_{S j}(x_i^{S}))_{i=1,\;j=1}^{P_{S},\;J_{S}}\\
\end{pmatrix},\\
\\
\end{aligned}
\end{equation*}
\begin{equation*}
\begin{aligned}
&A_2=\begin{pmatrix}
(\mathcal{B}\phi_{1 j}(x_0))_{j=1}^{J_1} & \ldots& 0\\
0 & \ldots & (\mathcal{B}\phi_{S j}(x_1))_{j=1}^{J_S}
\end{pmatrix},\\
\\
&A_3=\begin{pmatrix}
-(\phi_{1 j}(x_{P_1}^{1}))_{j=1}^{J_1} & (\phi_{2 j}(x_{1}^{2}))_{j=1}^{J_2} & 0 & \ldots & 0\\
0 & -(\phi_{2 j}(x_{P_2}^{2}))_{j=1}^{J_2} & (\phi_{3 j}(x_{1}^{3}))_{j=1}^{J_3} & \ldots & 0\\
\vdots & \vdots & \ddots & \ddots & \vdots\\
0 & 0 & 0 & -(\phi_{S-1\ j}(x_{P_{S-1}}^{S-1}))_{j=1}^{J_{S-1}} & (\phi_{S j}(x_{1}^{S}))_{j=1}^{J_S}\\
\end{pmatrix},\\
\\
\end{aligned}
\end{equation*}
\begin{equation*}
\begin{aligned}
&A_4=\begin{pmatrix}
-(\frac{d\phi_{1 j}}{d x}(x_{P_1}^{1}))_{j=1}^{J_1} & (\frac{d\phi_{2 j}}{d x}(x_{1}^{2}))_{j=1}^{J_2} & 0 & \ldots & 0\\
0 & -(\frac{d\phi_{2 j}}{d x}(x_{P_2}^{2}))_{j=1}^{J_2} & (\frac{d\phi_{3 j}}{d x}(x_{1}^{3}))_{j=1}^{J_3} & \ldots & 0\\
\vdots & \vdots & \vdots & \ddots & \vdots\\
0 & 0 & 0 & -(\frac{d\phi_{S-1\ j}}{d x}(x_{P_{S-1}}^{S-1}))_{j=1}^{J_{S-1}} & (\frac{d\phi_{S j}}{d x}(x_{1}^{S}))_{j=1}^{J_S}\\
\end{pmatrix},\\
\\
&\boldsymbol{v_1}=\begin{pmatrix}
f(x_1^1)\\
\ldots\\
f(x_{J_1}^1)\\
\ldots\\
\ldots\\
f(x_{J_S}^S)\\
\end{pmatrix},\quad\quad\quad
\boldsymbol{v_2}=\begin{pmatrix}
g(x_0)\\
g(x_1)\\
\end{pmatrix},\quad\quad\quad
\boldsymbol{v_3}=\begin{pmatrix}
0\\
\ldots\\
0\\
\end{pmatrix},\quad\quad\quad
\boldsymbol{v_4}=\begin{pmatrix}
0\\
\ldots\\
0\\
\end{pmatrix}\,.
\end{aligned}
\end{equation*}

Here are some elucidations:
\begin{itemize}
\item{$A_1\in \mathbb{R}^{P\times M}$ and $\boldsymbol{v_1}\in \mathbb{R}^{P}$ correspondingly incorporate the condition for satisfying the strong form of equation}
\item{$A_2\in \mathbb{R}^{2\times M}$ and $\boldsymbol{v_2} \in \mathbb{R}^{2}$ impose the boundary condition on the function}
\item{$A_3\in \mathbb{R}^{(S-1)\times M}$ and $\boldsymbol{v_3}\in \mathbb{R}^{S-1}$ impose continuity conditions at boundaries of blocks}
\item{Finally, $A_4\in \mathbb{R}^{(S-1)\times M}$ and $v_4\in \mathbb{R}^{S-1}$ put differentiability conditions at boundaries of blocks}
\end{itemize}
Hereby, problem of finding a strong solution of the equation \eqref{eq:1} was reformulated to linear system \eqref{eq:2}, which can be solved by various conventional methods.

\section{Weak RFM}
In this section, we provide a detailed introduction to the main ideas and implementation details of WRFM, highlighting its differences from the standard RFM.
\label{sec:new}
\subsection{Weak Formulation}
\label{subsec:weak_formulation}
\label{weak_formulation}
\begin{equation}\label{eq:3}
\left\{
\begin{aligned}
&\int_{\Omega} u(\boldsymbol{x})\cdot\mathcal{L}\varphi(\boldsymbol{x})d\boldsymbol{x}= \int_{\Omega}f(\boldsymbol{x})\cdot\varphi(\boldsymbol{x})d\boldsymbol{x}\,,\\
&\int_{\partial\Omega} u(\boldsymbol{x})\cdot\mathcal{B}\tilde{\varphi}(\boldsymbol{x})d\boldsymbol{x}= \int_{\partial\Omega}g(\boldsymbol{x})\cdot\tilde{\varphi}(\boldsymbol{x})d\boldsymbol{x}\,.\\
\end{aligned}
\right.
\end{equation}
According to the most common definition of weak solutions, $u(x)$ is weak solution of \eqref{eq:1} if equation \eqref{eq:3} holds for all $\varphi \in \mathcal{D}(\mathbb{R}^d)$ with $\operatorname{supp}(\varphi) \subset \Omega$, and for all $\tilde{\varphi} \in \mathcal{D}(\mathbb{R}^{d-1})$ with $\operatorname{supp}(\tilde{\varphi}) \subset \partial\Omega$, where:

\begin{equation*}
\begin{aligned}
&\mathcal{L}\varphi(\boldsymbol{x}) = \sum_{\boldsymbol{\alpha}\in\mathbb{N}^{d}} (-1)^{|\boldsymbol{\alpha}|}\cdot\partial^{\boldsymbol{\alpha}}(a_{\boldsymbol{\alpha}}(\boldsymbol{x})\varphi(\boldsymbol{x}))\,,\\
&\mathcal{B}\tilde{\varphi}(\boldsymbol{x}) = \sum_{\boldsymbol{\beta}\in\mathbb{N}^{d}} (-1)^{|\boldsymbol{\beta}|}\cdot\partial^{\boldsymbol{\beta}}(b_{\boldsymbol{\beta}}(\boldsymbol{x})\tilde{\varphi}(\boldsymbol{x}))\,.
\end{aligned}
\end{equation*}
It is noticeable that in formulation \eqref{eq:3} there are no requirements regarding the differentiability of solution $u(x)$, therefore space of weak solutions is larger than the original one. However, it is practically troublesome to prove that function is a weak solution by definition, since \eqref{eq:3} should be proven for infinite set of test functions. In order to overcome this difficulty, well known Weierstrass' theorem on the approximation of functions by trigonometric polynomials \cite{weierstrass1885} can be used. According to this theorem, if $\varphi\in C_{[-\pi,\pi]}$ and $\varphi(-\pi)=\varphi(\pi)$, then $\forall \varepsilon > 0,\exists N \in \mathbb{N}:$
\begin{equation*}
\abs{\varphi(x)-\left(\sum_{k=1}^{N}a_k \sin(kx) + b_k \cos(kx)\right)}<\varepsilon,\quad\forall x\in[-\pi,\pi]\,.
\end{equation*}
If $\varphi\in C_{[0,\pi]}$ and $\varphi(0)=\varphi(\pi)=0$, then we can expand it to odd $\tilde\varphi\in C_{[-\pi,\pi]}$ and apply the same theorem:
\begin{equation*}
\abs{\tilde\varphi(x)-\left(\sum_{k=1}^{N}\tilde a_k \sin(kx) + \tilde b_k \cos(kx)\right)}<\varepsilon,\quad\forall x\in[-\pi,\pi]\,.
\end{equation*}
Since $\tilde{\varphi}(x)$ is odd, $\cos(kx)$ is even, and $\tilde{b}_k = 0$ for all $k$. Hence:
\begin{equation*}
\abs{\varphi(x)-\sum_{k=1}^{N}\tilde a_k \sin(kx)}<\varepsilon,\quad\forall x\in[0,\pi]\,.
\end{equation*}
Consequently, for an arbitrary domain $[x_0, x_1]$ and the function $\varphi(x) \in C([x_0, x_1])$ with $\varphi(x_0) = \varphi(x_1) = 0$, the following expression is justified:

\begin{equation*}
\varphi(x) = \sum_{k=1}^{+\infty} c_k \sin\left(\frac{\pi k (x-x_0)}{x_1-x_0}\right)\,.
\end{equation*}
In the context of our problem, this implies that representation \eqref{repr:1} is valid for all test functions $\varphi \in \mathcal{D}(\mathbb{R})$ such that $\operatorname{supp}(\varphi) \subset [x_0, x_1]$:
\begin{equation}\label{repr:1}
\varphi(x) = w_{[x_0,x_1]}\cdot\sum_{k=1}^{+\infty} c_k \sin\left(\frac{\pi k (x-x_0)}{x_1-x_0}\right)\,,
\end{equation}
where $w_{[x_0,x_1]}$ is a window function ensuring that $supp(\varphi(x))\subset[x_0,x_1]$. Theoretically $supp(w_{[x_0,x_1]})=(x_0,x_1)$, but in numerical examples $\phi_{[x_0,x_1]}^{c}(x)$ can be used.\\
In 2d case, similar expression can be obtained considering common representation of multidimensional functions:
$$\varphi(x,y) = \sum_{n} \large{X}_n(x)\cdot \large{Y}_n(y)\,.$$
Using \eqref{repr:1} for both $\large{X}_n(x)$ and $\large{Y}_n(y)$:
\begin{equation}\label{repr:2}
\varphi(x,y) = w_{[x_0,x_1]\times[y_0,y_1]}\cdot\sum_{k,l} \tilde c_{kl} \sin\left(\frac{\pi k (x-x_0)}{x_1-x_0}\right)\cdot \sin\left(\frac{\pi l (y-y_0)}{y_1-y_0}\right)\,.
\end{equation}
Therefore, instead of proving \eqref{eq:3} for all $\varphi(\boldsymbol{x}) \in \mathcal{D}(\mathbb{R}^{d})$ with $\operatorname{supp}(\varphi) \subset \Omega$, it suffices to verify it for a sufficiently rich set of functions $\{\varphi_k\}$ or $\{\varphi_{kl}\}$ in the one- and two-dimensional cases, respectively:
\begin{equation*}
\begin{aligned}
\varphi_k&=w_{[x_0,x_1]}\cdot\sin\left(\frac{\pi k (x-x_0)}{x_1-x_0}\right)\,,\\
\varphi_{kl}&=w_{[x_0,x_1]\times[y_0,y_1]}\cdot\sin\left(\frac{\pi k (x-x_0)}{x_1-x_0}\right)\cdot \sin\left(\frac{\pi l (y-y_0)}{y_1-y_0}\right)\,.
\end{aligned}
\end{equation*}
For higher dimensions similar expressions can also be derived.

\subsection{Loss Function and Optimization for Weak Formulation}
\label{subsec:weak_rfm}
This section presents loss function and optimization process for the case of weak formulation defined in \eqref{eq:3}, introducing the algorithm of Weak RFM. According to weak form of equation, loss function has the form:
\begin{equation*}
\begin{aligned}
    \text { Loss }=\sum_{\boldsymbol{k} \in \mathbb{N}^d}^{|k|\leq K_I} \sum_{i=1}^{E_I} \left\| \int_{\Omega} u(\boldsymbol{x})\cdot\mathcal{L}^{i}\varphi_{\boldsymbol{k}}(\boldsymbol{x})- f(\boldsymbol{x})^{i}\cdot\varphi_{\boldsymbol{k}} (\boldsymbol{x})d\boldsymbol{x}\right\|_{l^2}^2 \\
    +\sum_{\boldsymbol{n} \in \mathbb{N}^{d-1}}^{|n|\leq K_B} \sum_{l=1}^{E_B} \left\| \int_{\partial\Omega} u(\boldsymbol{x})\cdot\mathcal{B}^{l}\tilde{\varphi}_{\boldsymbol{n}}(\boldsymbol{x})- g(\boldsymbol{x})^{l}\cdot\tilde{\varphi_{\boldsymbol{n}} }(\boldsymbol{x}) d\boldsymbol{x}\right\|_{l^2}^2\,.
    \end{aligned}
\end{equation*}
And exactly as in case of strong formulation, this optimization problem can be reformulated as a linear system and solved by large spectrum of  algorithms. In case $d=1$, $E_I=1$, $E_B=1$ and $\mathcal{B} u=u$ the system has the form \eqref{eq:4}.

\begin{equation}\label{eq:4}
\mathcal{A} \times \boldsymbol{u}=\boldsymbol{w}\,,
\end{equation}
where
\begin{equation*}
\begin{aligned}
\mathcal{A}=\begin{pmatrix}
\mathcal{A}_{1}\\
\mathcal{A}_{2}\\
\mathcal{A}_{3}
\end{pmatrix},&\quad
\boldsymbol{u}=\begin{pmatrix}
u_{11}\\
u_{21}\\
\vdots \\
u_{J_{S} S}
\end{pmatrix},\quad
\boldsymbol{w}=\begin{pmatrix}
\boldsymbol{w_{1}} \\
\boldsymbol{w_{2}} \\
\boldsymbol{w_{3}}
\end{pmatrix}\,.
\end{aligned}
\end{equation*}

We denote amount of test functions $\varphi^n$ corresponding to the $n$-th block as $K_n$ and total amount of test functions as $K=\sum_{n=1}^S K_n$. At the same time, there is no need in test functions for boundary conditions $\tilde \varphi$ since $\mathcal{B}$ does not contain any derivatives in the considered case. Therefore, structure of $\mathcal{A}$ and $\boldsymbol{w}$ takes the form as follows:

\begin{equation*}
\begin{aligned}
&\mathcal{A}_1=\begin{pmatrix}
(\int_{\Omega_1} \phi_{1 j}(x)\cdot\mathcal{L}\varphi_{i}^{1}(x)dx)_{i=1,\;j=1}^{K_1,\;J_1} & \ldots & 0\\
\vdots & \ddots & \vdots\\
0 & \ldots & (\int_{\Omega_S} \phi_{S j}(x)\cdot\mathcal{L}\varphi_{i}^{S}(x)dx)_{i=1,\;j=1}^{K_S,\;J_S}\\
\end{pmatrix},\\
\\
&\mathcal{A}_2=\begin{pmatrix}
(\phi_{1 j}(x_0))_{j=1}^{J_1} & \ldots& 0\\
0 & \ldots & (\phi_{S j}(x_1))_{j=1}^{J_S}
\end{pmatrix},\\
\\
\end{aligned}
\end{equation*}
\begin{equation*}
\begin{aligned}
&\mathcal{A}_3=\begin{pmatrix}
-(\phi_{1 j}(x_{P_1}^{1}))_{j=1}^{J_1} & (\phi_{2 j}(x_{1}^{2}))_{j=1}^{J_2} & 0 & \ldots & 0\\
0 & -(\phi_{2 j}(x_{P_2}^{2}))_{j=1}^{J_2} & (\phi_{3 j}(x_{1}^{3}))_{j=1}^{J_3} & \ldots & 0\\
\vdots & \vdots & \ddots & \ddots & \vdots\\
0 & 0 & 0 & -(\phi_{S-1\ j}(x_{P_{S-1}}^{S-1}))_{j=1}^{J_{S-1}} & (\phi_{S j}(x_{1}^{S}))_{j=1}^{J_S}\\
\end{pmatrix},\\
\\
\end{aligned}
\end{equation*}
\begin{equation*}
\begin{aligned}
&\boldsymbol{w_1}=\begin{pmatrix}
\int_{\Omega_1} f(x)\cdot\varphi_{1}^{1}(x)dx\\
\ldots\\
\int_{\Omega_1} f(x)\cdot\varphi_{K_1}^{1}(x)dx\\
\ldots\\
\ldots\\
\int_{\Omega_S} f(x)\cdot\varphi_{K_S}^{S}(x)dx\\
\end{pmatrix},\quad\quad\quad\quad
\boldsymbol{w_2}=\begin{pmatrix}
g(x_0)\\
g(x_1)\\
\end{pmatrix},\quad\quad\quad\quad
\boldsymbol{w_3}=\begin{pmatrix}
0\\
\ldots\\
0\\
\end{pmatrix}.
\end{aligned}
\end{equation*}
It is worth mentioning, that $\varphi_k^n(x)$ is defined exactly as in expression \eqref{repr:1}:
$$\varphi_k^n=w_{[x_0^n,x_1^n]}\cdot\sin\left(\frac{\pi k (x-x_0^n)}{x_1^n-x_0^n}\right)\,,$$ 
where $x_0^n$ and $x_1^n$ are boundaries of $\Omega_n$. In comparison with the standard RFM, dimensionality of matrices is different:
\begin{itemize}
\item{$\mathcal{A}_1\in \mathbb{R}^{K\times M}$ and $\boldsymbol{w_1}\in \mathbb{R}^{K}$ correspondingly impose the condition for satisfying the weak form of equation}
\item{$\mathcal{A}_2\in \mathbb{R}^{2\times M}$ and $\boldsymbol{w_2}\in \mathbb{R}^{2}$ impose the boundary condition on the function}
\item{$\mathcal{A}_3\in \mathbb{R}^{(S-1)\times M}$ and $\boldsymbol{w_3}\in \mathbb{R}^{S-1}$ put continuity conditions at boundaries of blocks, it is also noteworthy that there are no differentiability conditions}
\end{itemize}
This algorithm was implemented for up to 3D cases for both simple and complex geometry. Results are presented in the next section.

\section{Numerical Results}\label{sec:num_results}
This section evaluates the performance of Weak RFM by comparing it to two popular methods, namely PINN and WAN. PINN inherently has the ability to find weak solutions, as reformulating the equation into an optimization problem enables it to identify solutions that satisfy the equation almost everywhere—a concept closely aligned with the definition of weak solutions. WAN, introduced in \cite{zang2020wan}, builds upon the PINN framework but incorporates the weak formulation of the equation and employs an adversarial neural network as the test function. It is reported to be more effective than PINN in finding weak solutions, particularly in higher-dimensional cases.
To assess the performance of these methods, we compare their accuracy, consuming time, and the number of  parameters. Specifically, $L_2$ and $L_{\infty}$ norms are used to quantify accuracy:
\begin{equation*}
L_2^u:=\frac{\left\|u-u_e\right\|_2}{\left\|u_e\right\|_2}, \quad L_{\infty}^u:=\frac{\left\|u-u_e\right\|_{\infty}}{\left\|u_e\right\|_{\infty}},
\end{equation*}
where $u_e$ is the reference solution.
\subsection{Experimental Setup}
The PyTorch implementation of PINN is employed, incorporating standard enhancements such as exponential learning rate decay and the use of the L-BFGS optimizer following initial training with ADAM. For fair comparison, WAN is also reimplemented in PyTorch. All neural network training is performed on a GPU (Tesla V100-PCIE-32GB), whereas WRFM solves its linear system exclusively on a CPU (Intel Xeon CPU E5-2630 v4 @ 2.20GHz). Optimal hyperparameters for each method are carefully selected to ensure a balance between accuracy and comparable training time. The best-performing hyperparameters for WRFM are summarized in Table \ref{tab:hyperparameters}. Each experiment is repeated three times, and the reported results represent the averaged outcomes unless otherwise stated.

\subsection{2D Helmholtz equation}\label{Helmholtz}
Helmholtz equation is considered on regular domain $\Omega=\{(x,y)\in\mathbb{R}^2|\ x\in[-1,1], y\in[0,1]\}$ for $\lambda=1$:
\begin{equation}\label{eq:Helmholtz}
\left\{
\begin{aligned}
&\Delta u + \lambda u = \sinh(|x|)((1+\lambda-4y^2)\cos(y^2)-2\sin(y^2))\,,\\
&u|_{(x,y)\in\partial\Omega} = \sinh(|x|)\cos(y^2)\,.
\end{aligned}
\right.
\end{equation}

\begin{figure}[h]
\centering
\includegraphics[width=3.7in]{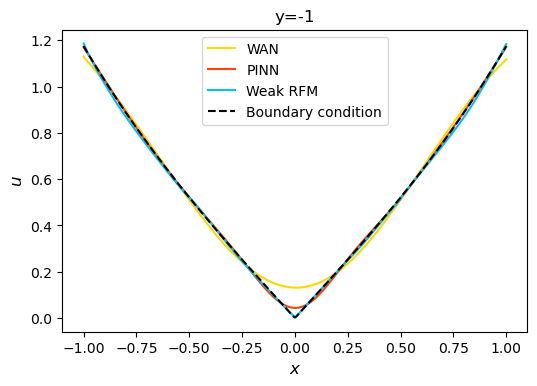}
\caption{Helmholtz equation: comparison of obtained solutions on boundary.}
\label{pic:Helmholtz_bc}
\end{figure}

\begin{figure}[H]
\centering
\begin{subfigure}{.5\textwidth}
  \centering
  \includegraphics[width=2.3in]{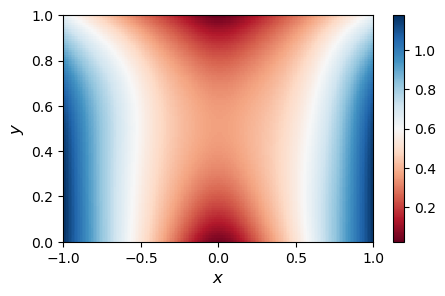}
  \caption{PINN solution}
\end{subfigure}%
\begin{subfigure}{.5\textwidth}
  \centering
  \includegraphics[width=2.3in]{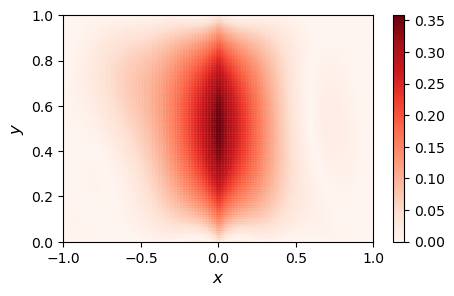}
  \caption{PINN absolute error}
\end{subfigure}
\begin{subfigure}{.5\textwidth}
  \centering
  \includegraphics[width=2.3in]{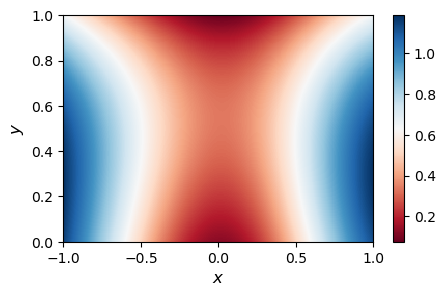}
  \caption{WAN solution}
\end{subfigure}%
\begin{subfigure}{.5\textwidth}
  \centering
  \includegraphics[width=2.3in]{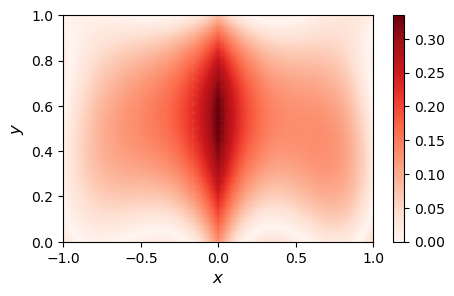}
  \caption{WAN absolute error}
  \label{fig:sub4}
\end{subfigure}
\begin{subfigure}{.5\textwidth}
  \centering
  \includegraphics[width=2.3in]{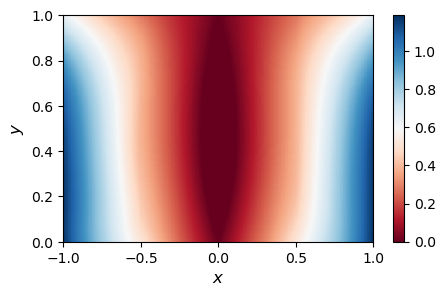}
  \caption{WRFM solution}
\end{subfigure}%
\begin{subfigure}{.5\textwidth}
  \centering
  \includegraphics[width=2.3in]{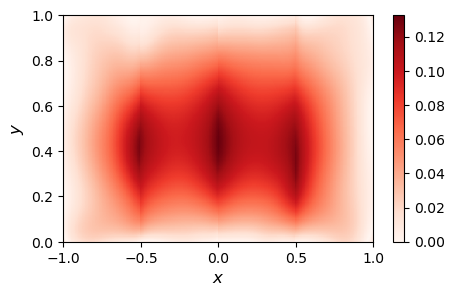}
  \caption{WRFM absolute error}
  \label{fig:sub4}
\end{subfigure}
\caption{Helmholtz equation on regular domain.}
\label{pic:Helmholtz_reg}
\end{figure}

In this case, there is no strong solution due to the lack of differentiability at $x=0$, but the equation still admits a weak solution. As illustrated in Fig. \ref{pic:Helmholtz_reg} and summarized in Table \ref{tab:helmholtz_reg}, we applied PINN, WAN, and WRFM to solve this problem, evaluating their performance using an identical number of parameters. In this example, WAN outperforms PINN by achieving comparable accuracy while reducing computational time. However, WRFM surpasses both methods, achieving superior accuracy in significantly less time and with fewer parameters. This advantage arises from WRFM's enhanced ability to handle the non-smooth features of the boundary conditions, as shown in Fig. \ref{pic:Helmholtz_bc}.\\

\begin{table}[H]
\centering
\begin{tabular}{|l|cccc|}
\hline
     \textbf{Method} & \textbf{Params} & \textbf{Time (s)} & \textbf{$L_2$} & \textbf{$L_\infty$} \\ \specialrule{.15em}{.0em}{.0em}
    PINN & 360 & 137 & 0,169 & 0,388 \\
        ~ & \textbf{500} & \textbf{106} & \textbf{0,157} & \textbf{0,359} \\
        ~ & 680 & 134 & 0,172 & 0,385 \\ \hline
    WAN & 360 & 57 & 0,213 & 0,429 \\
    ~ & \textbf{500} & \textbf{55} & \textbf{0,177} & \textbf{0,349} \\
    ~ & 680 & 54 & 0,209 & 0,403 \\ \hline
    WRFM & 80 & 17 & 0,066 & 0,146 \\
    ~ & 120 & 16 & 0,069 & 0,165 \\
    ~ & \textbf{160} & \textbf{15} & \textbf{0,091} & \textbf{0,133} \\ \hline
\end{tabular}
\caption{Performance of methods in solving the Helmholtz equation on regular domain.}
\label{tab:helmholtz_reg}
\end{table}

Equation \eqref{eq:Helmholtz} is also examined on a complex domain:

\begin{equation*}
\begin{aligned}
\Omega = \{(x,y)\in\mathbb{R}^2|\ x\in[-2,1], y\in[0,1]\} \setminus \bigcup_{i=1}^N D_r(x_i,y_i)\,,\\
\end{aligned}
\end{equation*}
where \( D_r(x_i, y_i) \) denotes a disk with radius \( r \) and center \( (x_i, y_i) \). In the case under consideration, \( r = 0.05 \), \( N = 8 \), and the disks are arranged in a hexagonal packing pattern at the center of the domain.

\begin{figure}[H]
\centering
\begin{subfigure}{.5\textwidth}
  \centering
  \includegraphics[width=2.3in]{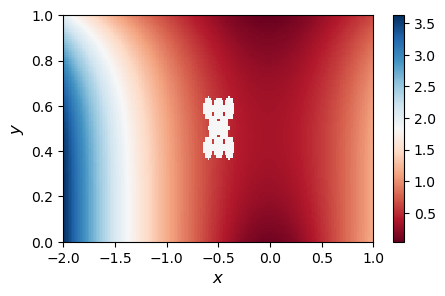}
  \caption{PINN solution}
\end{subfigure}%
\begin{subfigure}{.5\textwidth}
  \centering
  \includegraphics[width=2.3in]{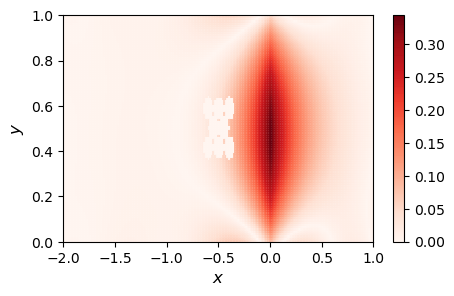}
  \caption{PINN absolute error}
\end{subfigure}
\begin{subfigure}{.5\textwidth}
  \centering
  \includegraphics[width=2.3in]{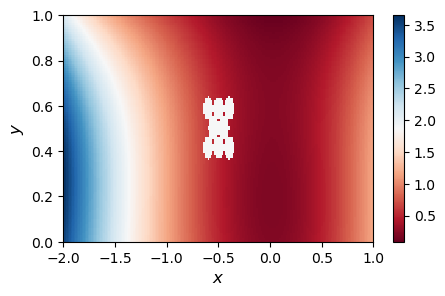}
  \caption{WAN solution}
\end{subfigure}%
\begin{subfigure}{.5\textwidth}
  \centering
  \includegraphics[width=2.3in]{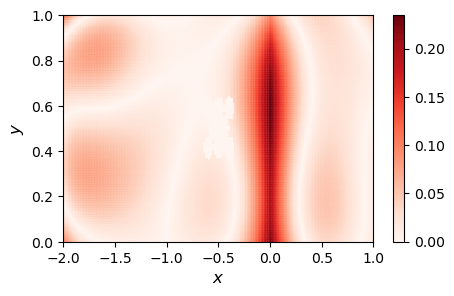}
  \caption{WAN absolute error}
\end{subfigure}
\begin{subfigure}{.5\textwidth}
  \centering
  \includegraphics[width=2.3in]{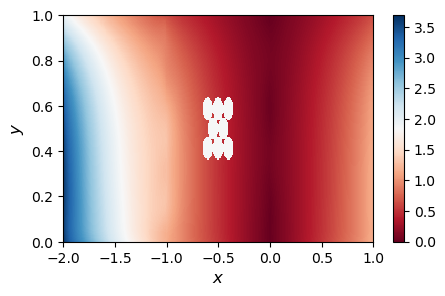}
  \caption{WRFM solution}
\end{subfigure}%
\begin{subfigure}{.5\textwidth}
  \centering
  \includegraphics[width=2.3in]{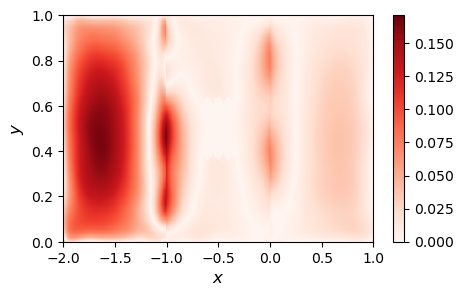}
  \caption{WRFM absolute error}
\end{subfigure}
\caption{Helmholtz equation on complex domain.}
\label{pic:Helmholtz_com}
\end{figure}

The considered methods are still capable of providing precise solutions, as shown in Figure \ref{pic:Helmholtz_com}, with WAN again outperforming PINN in terms of accuracy. However, the best result in this case is obtained by WRFM, which produces a more accurate solution using fewer trainable parameters and requiring fewer resources for training, as summarized in Table \ref{tab:helmholtz_com}. It is important to note, however, that the accuracy of WRFM tends to be influenced by the location, shape, and size of the cut-out area, whereas the other methods offer more stable but less accurate results.

\begin{table}[H]
\centering
\begin{tabular}{|l|cccc|}
\hline
     \textbf{Method} & \textbf{Params} & \textbf{Time (s)} & \textbf{$L_2$} & \textbf{$L_\infty$} \\ \specialrule{.15em}{.0em}{.0em}
        PINN & 540 & 172 & 0,143 & 0,344 \\
        ~ & 920 & 181 & 0,144 & 0,344 \\
        ~ & \textbf{1400} & \textbf{177} & \textbf{0,142} & \textbf{0,342} \\ \hline
        WAN & 540 & 175 & 0,201 & 0,381 \\
        ~ & \textbf{920} & \textbf{176} & \textbf{0,116} & \textbf{0,251} \\
        ~ & 1400 & 168 & 0,141 & 0,243 \\ \hline
        WRFM & 240 & 45 & 0,111 & 0,242 \\
        ~ & \textbf{300} & \textbf{46} & \textbf{0,101} & \textbf{0,171} \\
        ~ & 360 & 49 & 0,102 & 0,185 \\ \hline
\end{tabular}
\caption{Performance of methods in solving the Helmholtz equation on complex domain.}
\label{tab:helmholtz_com}
\end{table}

\subsection{2D Static heat equation}
Another example of a complex domain challenge is the following problem for the static heat equation on the complex domain \( \Omega = \{(x, y) \in \mathbb{R}^2 \mid x \in [0, 1], y \in [0, 1] \} \setminus \{(x, y) \in \mathbb{R}^2 \mid x \in \left[\frac{1}{3}, \frac{2}{3}\right], y \in \left[\frac{1}{3}, \frac{2}{3}\right] \} \):

\begin{equation}\label{eq:StaticHeat}
\left\{
\begin{aligned}
\Delta u =\ &9\cdot sgn((x-y)(y+x-1))\,,\\
u(x,y) &= \begin{cases} 
+1\ &\text{if}\ x\in\{0;1\}, y\in(0,1)\,,\\
-1\ &\text{if}\ x\in(0,1), y\in\{0;1\}\,,\\
0\ &\text{if}\ x\in\{\frac{1}{3};\frac{2}{3}\}, y\in(\frac{1}{3},\frac{2}{3})\,,\\
0\ &\text{if}\ x\in(\frac{1}{3},\frac{2}{3}), y\in\{\frac{1}{3};\frac{2}{3}\}\,.\\
\end{cases}\\
\end{aligned}
\right.
\end{equation}
This equation has no strong solution due to the incompatibility of the specified boundary conditions, which are not continuous. Nevertheless, a weak solution to this problem exists:
\begin{equation*}
u(x,y) = \left(\frac{9}{2}(x^2-x) + 1\right)\cdot\mathbb{I}_{[(x-y)(y+x-1)]} + \left(\frac{9}{2}(y-y^2) - 1\right)\cdot\mathbb{I}_{[(y-x)(y+x-1)]},
\end{equation*}
where
\begin{equation*}
\begin{aligned}
\ &\mathbb{I}_{[f(x,y)]} = \begin{cases}
1\ \text{if}\ f(x,y)>0\,,\\
0\ \text{if}\ f(x,y)\leq 0\,.
\end{cases}
\end{aligned}
\end{equation*}

\begin{figure}[H]
\centering
\begin{subfigure}{0.5\textwidth}
  \includegraphics[width=2.45in]{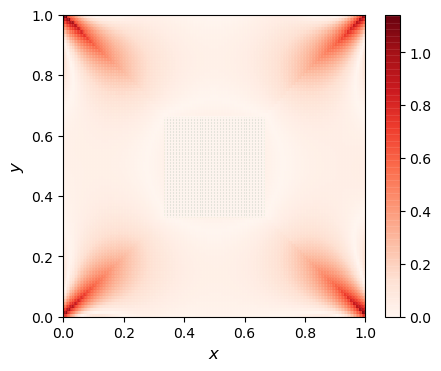}
  \caption{PINN}
\end{subfigure}%
\begin{subfigure}{0.5\textwidth}
  \includegraphics[width=2.45in]{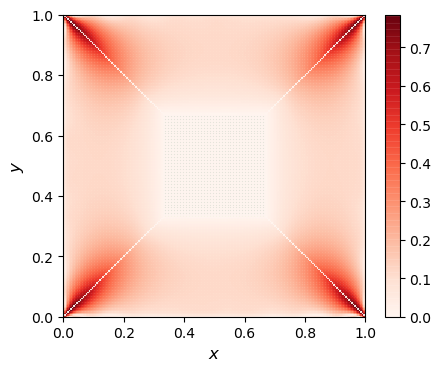}
  \caption{WRFM}
\end{subfigure}%
\caption{The absolute error for static heat equation on complex domain: (a) PINN, (b) WRFM.}
\label{pic:Static_Heat}
\end{figure}

Although WAN seems to fail in this case, the other two methods successfully solve the problem (see Figure \ref{pic:Static_Heat}) and demonstrate comparable accuracy. However, PINN utilized more computational resources to achieve similar results. It is notable that in instances where the analytical solution is discontinuous ($y=x$ or $y=1-x$), the resulting accuracy is diminished. The reason is that both methods are trying to find a continuous solution and approximate the break. This also helps to explain why $L_{\infty}$ metrics are high while $L_{2}$ are moderate in Table \ref{tab:static}.

\begin{table}[H]
\centering
\begin{tabular}{|l|cccc|}
\hline
    \textbf{Method} & \textbf{Params} & \textbf{Time (s)} & \textbf{$L_2$} & \textbf{$L_\infty$} \\ \specialrule{.15em}{.0em}{.0em}
    PINN & \textbf{920} & \textbf{382} & \textbf{0,164} & \textbf{1,189} \\
    ~ & 1980 & 388 & 0,164 & 1,190 \\
    ~ & 3440 & 389 & 0,164 & 1,205 \\ \hline
    WAN & 920 & 497 & 0,490 & 1,694 \\
    ~ & \textbf{1980} & \textbf{495} & \textbf{0,409} & \textbf{1,435} \\
    ~ & 3440 & 500 & 0,432 & 1,349 \\ \hline
    WRFM & 100 & 35 & 0,196 & 1,006 \\
    ~ & 300 & 54 & 0,184 & 0,807 \\
    ~ & \textbf{500} & \textbf{102} & \textbf{0,182} & \textbf{0,785} \\ \hline
\end{tabular}
\caption{Performance of methods in solving the static heat equation.}
\label{tab:static}
\end{table}

\subsection{3D Poisson`s equation}
\begin{figure}[H]
\centering
\begin{subfigure}{1.0\textwidth}
  \includegraphics[width=5.2in]{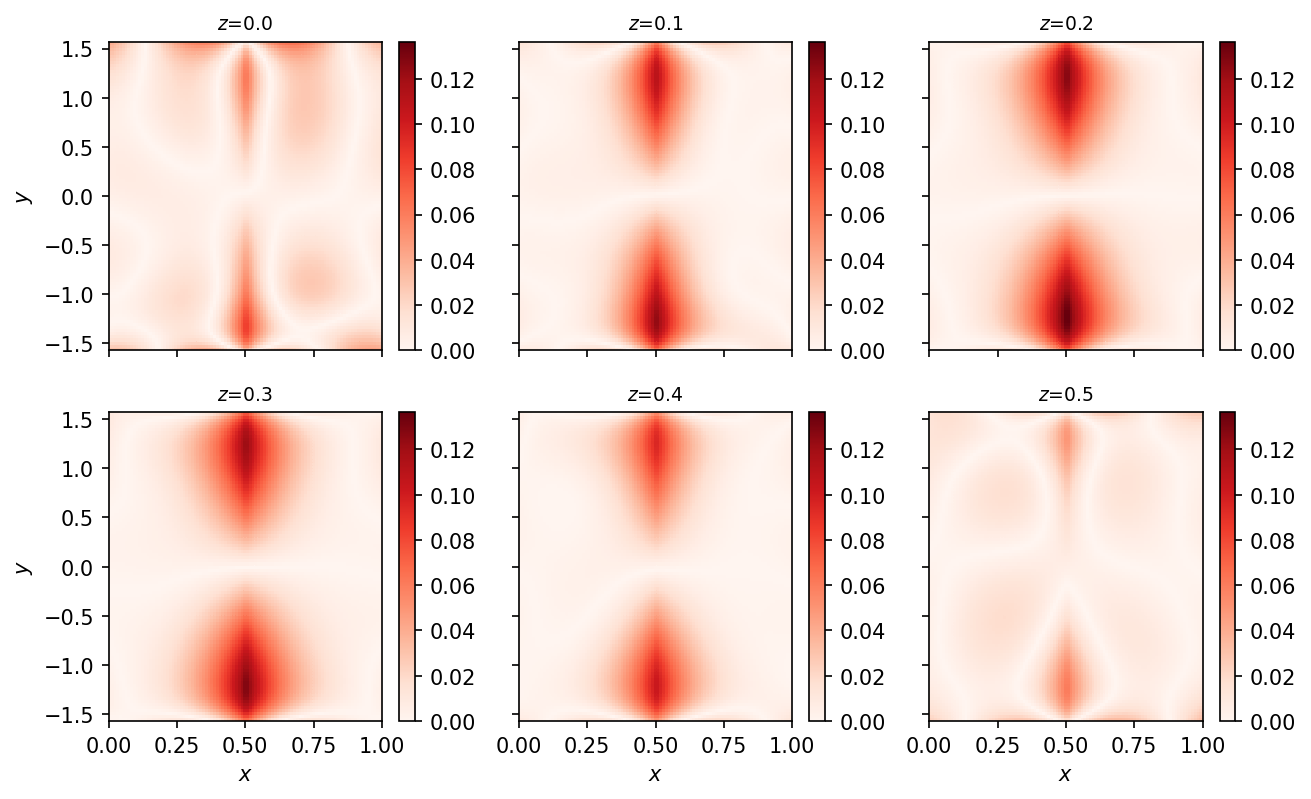}
  \caption{PINN}
\end{subfigure}
\begin{subfigure}{1.0\textwidth}
  \includegraphics[width=5.2in]{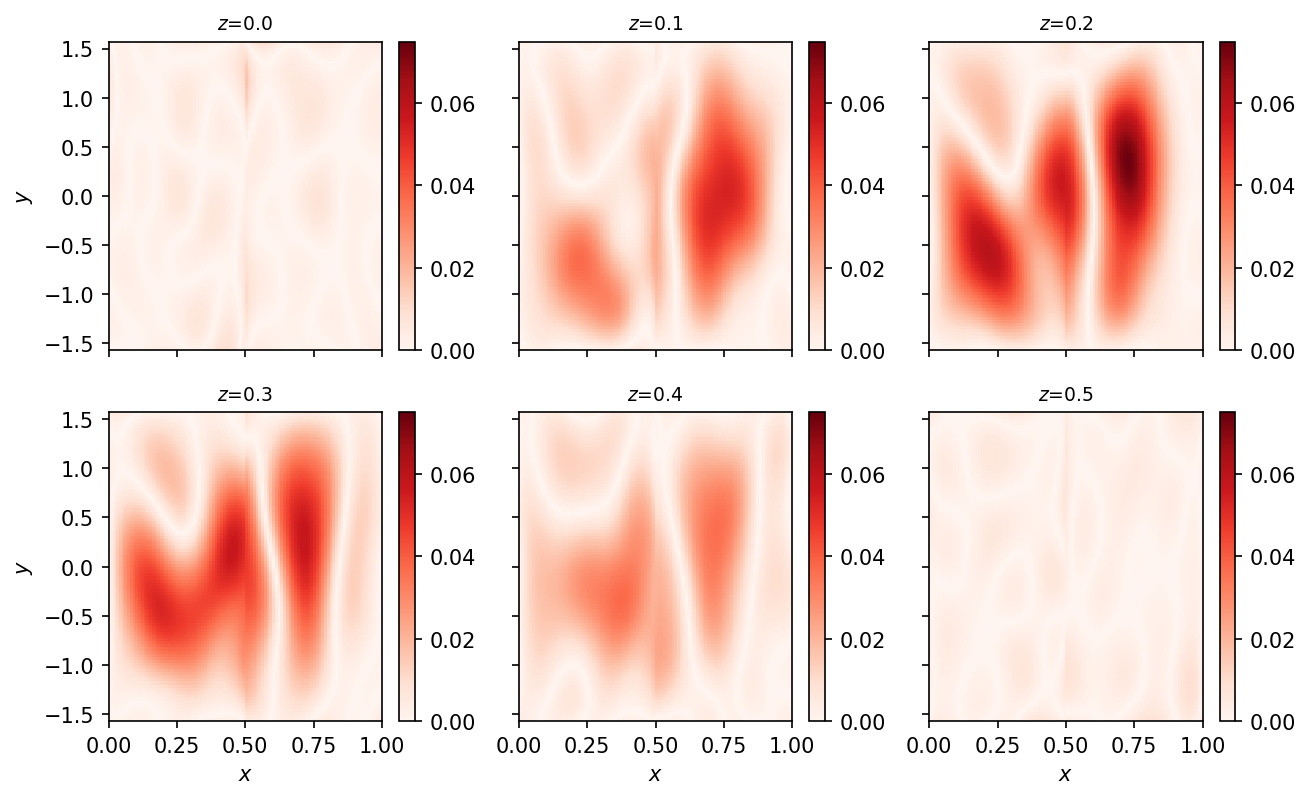}
  \caption{Weak RFM}
\end{subfigure}
\caption{The absolute error for Poisson`s equations of the two most accurate methods: (a) PINN, (b) WRFM.}
\label{pic:Poisson}
\end{figure}

Satisfying the Poisson equation does not necessarily imply the existence of a strong solution. Consider the following equation, in which the boundary condition is not smooth at $x=\frac{1}{2}$:
\begin{equation}\label{eq:Poisson}
\left\{
\begin{aligned}
&\Delta u = 2\sin(y)e^{-z}\,,\\
&u|_{(x,y,z)\in\partial\Omega} = \left(x-\mathbb{I}_{[x-\frac{1}{2}]}\right)^2 \sin(y) e^{-z}\,.
\end{aligned}
\right.
\end{equation}
Equation is considered on :
$$\Omega=\{(x,y,z)\in\mathbb{R}^3|\ x\in[0,1], y\in[-\frac{\pi}{2},+\frac{\pi}{2}], z\in[0,\frac{1}{2}]\}.$$
In 3D cases, the curse of dimensionality begins to affect Weak RFM. Given that a certain number of test functions must be selected for each dimension, the total size of the linear system increases, thus requiring more time for the calculation of \( A \) and \( w \), as well as for solving the linear system. Furthermore, since both steps rely on the CPU, this process takes a noticeable amount of time, comparable to the GPU-optimized PINN training process. Nevertheless, , as illustrated in Table \ref{tab:poisson}, in this case, WRFM finds a similar solution using fewer trainable parameters compared to both PINN and WAN. Moreover, it is notable that WRFM handles non-smooth regions much better than the other two methods , as illustrated in Figure \ref{pic:Poisson}.

\begin{table}[H]
\centering
\begin{tabular}{|l|cccc|}
\hline
    \textbf{Method} & \textbf{Params} & \textbf{Time (s)} & \textbf{$L_2$} & \textbf{$L_\infty$} \\ \specialrule{.15em}{.0em}{.0em}
    PINN & 940 & 510 & 0,038 & 0,136 \\
    ~ & 2010 & 491 & 0,037 & 0,136 \\
    ~ & \textbf{2940} & \textbf{616} & \textbf{0,036} & \textbf{0,133} \\ \hline
    WAN & 940 & 494 & 0,051 & 0,170 \\
    ~ & \textbf{2010} & \textbf{479} & \textbf{0,040} & \textbf{0,151} \\
    ~ & 2940 & 492 & 0,046 & 0,183 \\ \hline
    WRFM & 160 & 418 & 0,072 & 0,093 \\
    ~ & \textbf{200} & \textbf{438} & \textbf{0,062} & \textbf{0,074} \\
    ~ & 220 & 608 & 0,073 & 0,091 \\ \hline
\end{tabular}
\caption{Performance of methods in solving the Poisson`s equation.}
\label{tab:poisson}
\end{table}

\subsection{3D Heat equation}
\begin{figure}[H]
\centering
\begin{subfigure}{1.0\textwidth}
  \centering
  \includegraphics[width=5.2in]{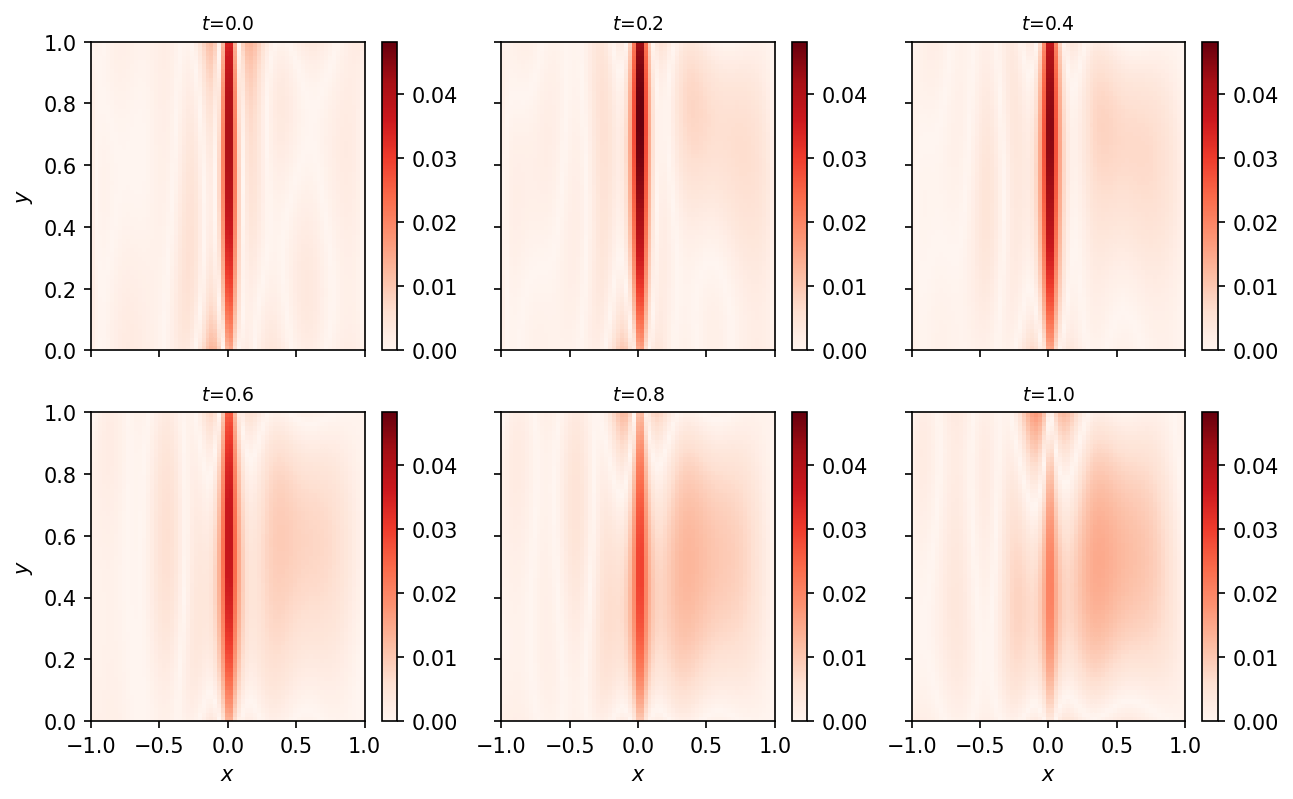}
  \caption{WAN}
\end{subfigure}
\begin{subfigure}{1.0\textwidth}
  \centering
  \includegraphics[width=5.2in]{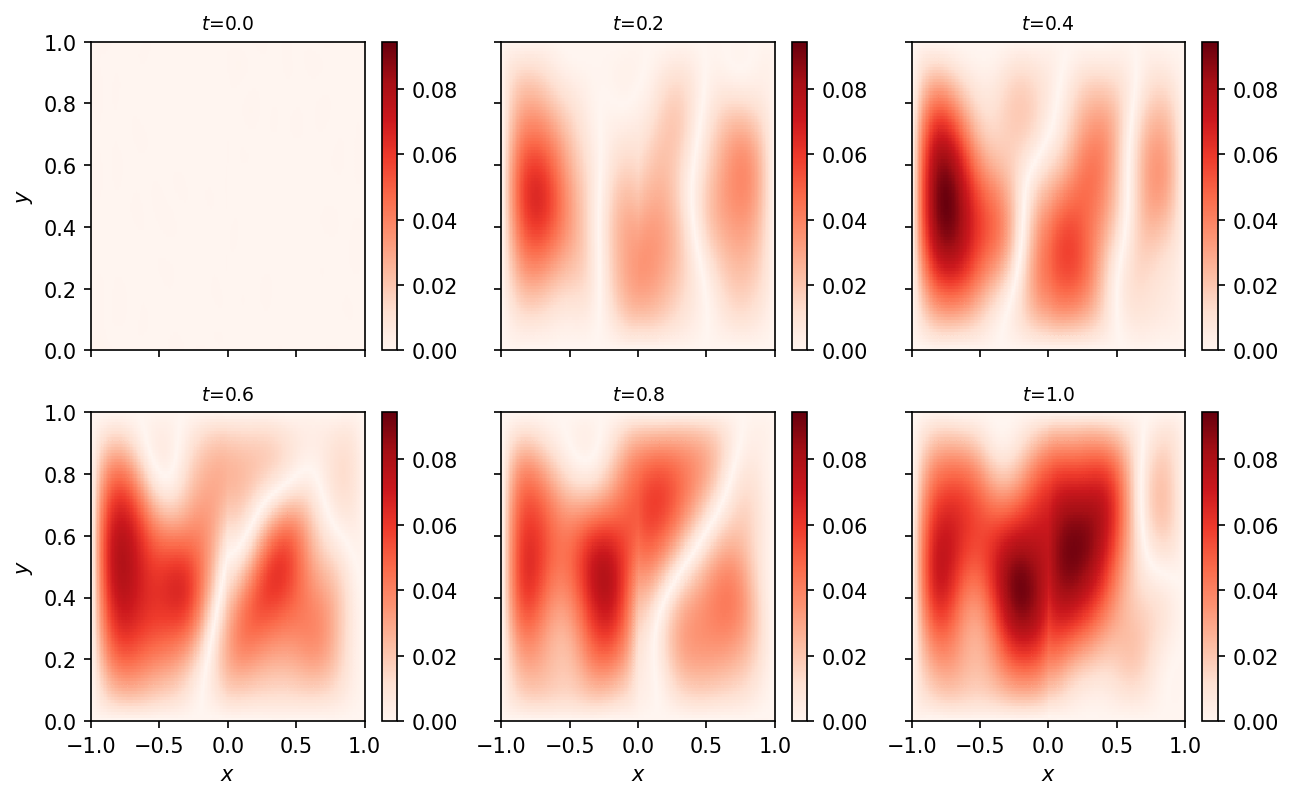}
  \caption{Weak RFM}
\end{subfigure}
\caption{The absolute error for Heat equation of the two most accurate methods: (a) PINN, (b) WRFM.}
\label{pic:Heat}
\end{figure}

Heat equation is considered on regular domain $\Omega=\{(x,y,t)\in\mathbb{R}^3|\ x\in[-1,1], y\in[0,1], t\in[0,1]\}$ for $a=1$:
\begin{equation}\label{eq:Heat}
\left\{
\begin{aligned}
&\Delta u - a^2 u_{t} = 2\sin\left(|x|+|y|-1\right)(t-1)e^{-\frac{t^2}{a^2}}\,,\\
&u|_{\partial\Omega\setminus\{(x,y,1)\}} = \sin\left(|x|+|y|-1\right)e^{-\frac{t^2}{a^2}}\,.
\end{aligned}
\right.
\end{equation}
In this case all methods take a considerable amount of time for finding solution. Performance of PINN is notably inferior while WRFM and WAN demonstrate significantly better results (see Table \ref{tab:heat}). It is noteworthy that in this particular case, the use of WAN appears to offer an advantage in terms of accuracy. Nevertheless, WRFM persists in demonstrating enhanced satisfaction with boundary conditions as well as superior handling of non-smooth parts of the solution (see Figure \ref{pic:Heat}).  

\begin{table}[H]
\centering
\begin{tabular}{|l|cccc|}
\hline
    \textbf{Method} & \textbf{Params} & \textbf{Time (s)} & \textbf{$L_2$} & \textbf{$L_\infty$} \\ \specialrule{.15em}{.0em}{.0em}
    PINN & 3480 & 1534 & 0,102 & 0,289 \\
    ~ & 5350 & 1488 & 0,102 & 0,289 \\
    ~ & \textbf{7900} & \textbf{2175} & \textbf{0,101} & \textbf{0,288} \\ \hline
    WAN & 3480 & 1251 & 0,016 & 0,066 \\
    ~ & \textbf{5350} & \textbf{1287} & \textbf{0,012} & \textbf{0,053} \\
    ~ & 7620 & 1253 & 0,016 & 0,052 \\ \hline
    WRFM & 200 & 867 & 0,045 & 0,129 \\
    ~ & 300 & 839 & 0,050 & 0,124 \\
    ~ & \textbf{400} & \textbf{1069} & \textbf{0,032} & \textbf{0,083} \\ \hline
\end{tabular}
\caption{Performance of methods in solving the heat equation.}
\label{tab:heat}
\end{table}

\begin{table}[H]
\centering
\begin{tabular}{|l|c|c|c|c|c|c|c|c|}
\hline
    \textbf{Problem} & $S$ & $J_n$ & $K_n^x$ & $K_n^y$ & $K_n^z$ & $P_n^x$ & $P_n^y$ & $P_n^z$ \\ \specialrule{.15em}{.0em}{.0em}
    2D Helmholtz, regular & 4 & 40 & 25 & 25 & - & 50 & 100 & - \\ \hline
    2D Helmholtz, complex & 3 & 100 & 25 & 15 & - & 150 & 100 & - \\ \hline
    2D static heat, complex & 1 & 500 & 50 & 50 & - & 100 & 100 & - \\ \hline
    3D Poisson`s, regular & 2 & 100 & 15 & 15 & 15 & 25 & 50 & 50 \\ \hline
    3D heat, regular & 2 & 200 & 15 & 15 & 15 & 50 & 50 & 50 \\ \hline
\end{tabular}
\caption{Optimal hyperparameters of Weak RFM in conducted experiments\\ ($S$ - amount of subdomains, $J_n$ - amount of feature functions in each subdomain,\\
$K_n = K_n^x \cdot K_n^y \cdot K_n^z$ - amount of test functions used in each subdomain,\\ 
$P_n^x,\ P_n^y,\ P_n^z$ - partitions of axes used for boundary conditions in each subdomain).}
\label{tab:hyperparameters}
\end{table}

\subsection{Sensitivity analysis}

\begin{figure}[H]
\centering
\begin{subfigure}{1.\textwidth}
  \centering
  \includegraphics[width=3.5in]{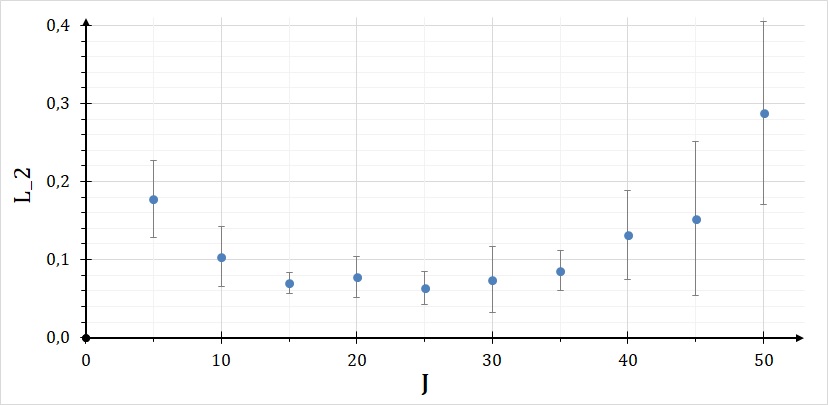}
\end{subfigure}
\begin{subfigure}{1.\textwidth}
  \centering
  \includegraphics[width=3.5in]{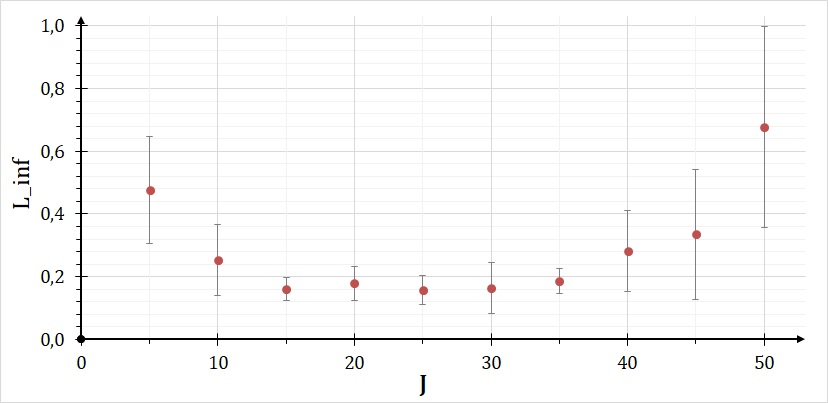}
\end{subfigure}
\caption{Dependence of accuracy on the number of basis functions.}
\label{pic:performance_1}
\end{figure}

\begin{figure}[H]
\centering
\begin{subfigure}{1.\textwidth}
  \centering
  \includegraphics[width=3.5in]{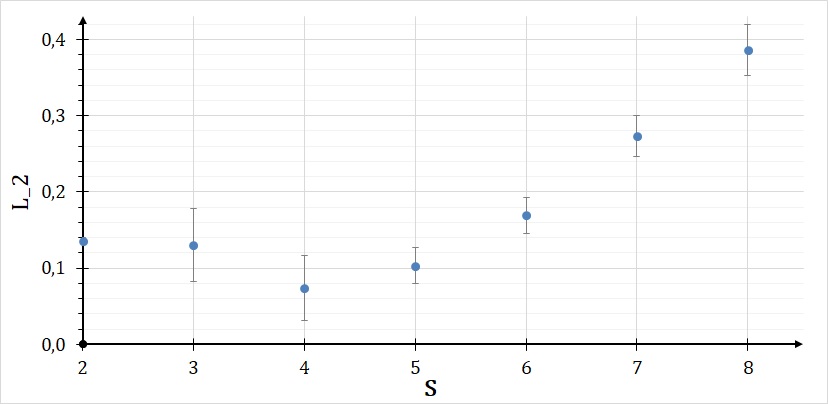}
\end{subfigure}
\begin{subfigure}{1.\textwidth}
  \centering
  \includegraphics[width=3.5in]{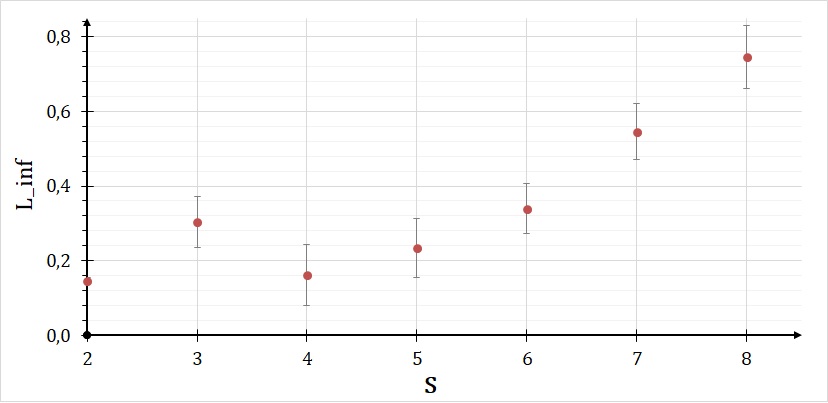}
\end{subfigure}
\caption{Dependence of accuracy on the number of domains.}
\label{pic:performance_2}
\end{figure}

The performance of the proposed method is highly dependent on its hyperparameter settings. While some parameters are less critical and can be fixed universally, two are particularly sensitive: $S$, the number of subdomains and $J$, the number of random feature functions used for each subdomain. Together, they determine the total number of unknowns $M=J\cdot S$ in linear system \eqref{eq:4}, significantly influencing accuracy. Figures \ref{pic:performance_1}, \ref{pic:performance_2} illustrate the relationships between the primary metrics $L_2,\ L_{\infty}$ and these parameters for the Helmholtz equation discussed in subsection \ref{Helmholtz}. All experiments were repeated 10 times, and the standard deviation was calculated to estimate variability.

It is evident that accuracy decreases significantly when $M=J\cdot S$ becomes excessively high or low. However, there exists a substantial range of $J$ values that yield high accuracy, making it relatively straightforward to identify the optimal configuration. While the accuracy of Weak RFM is sensitive to hyperparameter selection, this dependency is predictable and features a broad optimal range. As a result, this characteristic cannot be regarded as a serious limitation.

\label{sec:nr}
\section{Conclusion and Discussion}\label{sec:conclusion}
In this paper, we introduced and investigated Weak RFM, a novel method for finding weak solutions. Built upon the Random Feature Method framework, Weak RFM inherits all its advantages: it retains the mesh-free nature, leverages a least-squares formulation to solve linear systems, avoids the costly training process typical of neural networks, and maintains good accuracy. More importantly, Weak RFM extends the capability of RFM to problems involving weak solutions with low regularity—cases that standard RFM fails to handle.

We conducted a series of experiments on both 2D and 3D problems. The numerical results demonstrate the effectiveness of WRFM in handling PDEs lacking analytical solutions. Compared with existing ML-based approaches for finding weak solutions, such as PINN and WAN, WRFM achieves accurate results with significantly fewer parameters and reduced computational resources. This advantage is especially pronounced in 2D scenarios, where WRFM not only achieves higher accuracy but also runs substantially faster. It is worth noting that all WRFM experiments were performed using only CPUs, while the baseline ML methods utilized both CPUs and GPUs for training.

However, a current limitation of WRFM arises from its reliance on pre-defined test functions. As the dimensionality increases, the number of such functions grows exponentially, leading to scalability issues due to the resulting large linear systems. This phenomenon explains the increased computational time observed in 3D cases. Nevertheless, WRFM remains competitive and proves particularly effective in capturing the non-smooth regions of the solution.

In future work, we plan to explore the use of adversarially generated random feature functions as test functions. This approach may alleviate the curse of dimensionality and further enhance the performance of Weak RFM in high-dimensional settings.

\label{sec:con}
\section*{Conflict of interest statement}
There is no conflict of interest.

\section*{Acknowledgments}
The work is supported by National Key R\&D Program of China (No. 2022YFA1005200, No. 2022YFA1005202, and No. 2022YFA1005203), NSFC Major Research Plan -  Interpretable and General-purpose Next-generation Artificial Intelligence (No. 92270001 and No. 92270205), Anhui Center for Applied Mathematics, and the Major Project of Science \& Technology of Anhui Province (No. 202203a05020050).

\bibliographystyle{elsarticle-num} 
\bibliography{references}

\end{document}